\documentclass{amsart}
\usepackage[all]{xy}
\usepackage{latexsym,amssymb}
\numberwithin{equation}{section}
\newcommand\wtilde\widetilde
\newcommand\II{\mathrm{II}}
\newcommand\Vainerman{Va\u\i nerman}
\newcommand\Szlachanyi{Szlach\'anyi}
\newcommand{\E}{\operatorname E}
\newcommand{\lsum}{{\textstyle\sum}}
\newtheorem{Lem}{Lemma}[section]
\newtheorem{Prop}[Lem]{Proposition}
\newtheorem{Cor}[Lem]{Corollary}
\newtheorem{Thm}[Lem]{Theorem}

\theoremstyle{definition}
\newtheorem{Def}[Lem]{Definition}

\theoremstyle{remark}
\newtheorem{Rem}[Lem]{Remark}

\renewcommand\o{\otimes}

\newcommand\ol{\overline}
\newcommand\tR{\times_R}
\newcommand\tS{\times_S}
\newcommand\nt{\diamond}
\newcommand{\mapsfrom}{\mathrel{\mbox{$\leftarrow\joinrel\mapstochar$\,}}}
\DeclareMathOperator\Hom{\operatorname{Hom}}
\renewcommand\hom{\underline{\operatorname{hom}}}
\DeclareMathOperator\End{\operatorname{End}}
\DeclareMathOperator\id{\operatorname{id}}

\newcommand\op{{\operatorname{op}}}

\newcommand\C{\mathcal C}
\newcommand\bC{\mathbb C}

\newcommand\F{{\mathcal F}}

\newcommand\HMod[4]{{^{#1}_{#3}\mathcal M^{#2}_{#4}}}
\newcommand\LMod[1]{{_{#1}\mathcal M}}

\newcommand\Bimod[2]{{_{#1}\mathcal M_{#2}}}
\newcommand\BiMod[1]{\Bimod{#1}{#1}}
\newcommand\LComod[1]{{^{#1}\mathcal M}}

\newcommand\kmod{\mathcal M_k}

\newcommand{\leer}{\operatorname{--}}
\newcommand{\ou}[1]{\mathrel{\mathop{\otimes}_{#1}}}

\newcommand{\nti}[1]{\mathrel{\mathop{\diamond}_{#1}}}

\newcommand\sw[1]{{}_{(#1)}}
\newcommand\swm[1]{{}_{(-#1)}}

\newcommand\inv{^{-1}}
\renewcommand\epsilon\varepsilon

\newcommand\sqm{\ensuremath{\sqrt{\text{Morita}}}}
\newcommand\sqmsim{\overset{\sqrt{\text{M}}}\sim}
\newcommand\msim{\overset{\text{M}}\sim}
\makeatletter
\def\namelabel#1#2{\@bsphack
  \protected@write\@auxout{}%
         {\string\newlabel{#1.nme}{{#2}{#2}}}%
  \@esphack}
\def\nmlabel#1#2{\label{#2}\namelabel{#2}{#1}}
\newcommand\nmref[1]{\ref{#1.nme}\ \ref{#1}}
\makeatother

\newlength{\itemwidth}
\begin{document}
\title[Morita base change]%
  {Morita base change in quantum groupoids}
\author{Peter Schauenburg}
\address{Mathematisches Institut der Universit\"at M\"unchen, 
Theresienstr.~39, 80333~M\"unchen, Germany}
\email{schauen@rz.mathematik.uni-muenchen.de}
\begin{abstract}
  Let $L$ be a quantum semigroupoid, more precisely a $\tR$-bialgebra
  in the sense of Takeuchi. We describe a procedure replacing the
  algebra $R$ by any Morita equivalent, or in fact more generally
  any \sqm\ equivalent (in the sense of Takeuchi) algebra $S$ to
  obtain a $\tS$-bialgebra $\wtilde H$ with the same monoidal 
  representation category. 
\end{abstract}
\maketitle
\section{Introduction}
Quantum groupoids (or Hopf algebroids) are algebraic structures designed
to be the analogs of (the function algebras of) groupoids in the
realm of noncommutative geometry. A groupoid consists of a 
set $G$ of arrows, and a set $V$ of vertices. Thus a quantum groupoid
consists of an algebra $L$ (the function algebra on the noncommutative
space of arrows) and an algebra $R$ (the function algebra on the 
noncommutative space of vertices. The assignment to an arrow of its
source and target vertices defines two maps $G\rightrightarrows V$.
Thus the definition of a quantum groupoid involves two maps
$R\rightrightarrows L$; it turns out to be the right choice to assume
one of these to be an algebra, the other
an anti-algebra map, and to assume that the images of the two commute.
Since multiplication in $G$ is an only partially defined map, 
comultiplication in $L$ maps from $L$ to some tensor product
$L\ou RL$; one has to make the right choice of module structures to
define the tensor product, and one needs to assume that comultiplication
actually maps to a certain subspace of $L\ou RL$ to be able to state
that comultiplication is assumed to be an algebra map.

The first version of a quantum (semi)groupoid or bialgebroid or
Hopf algebroid was considered by Takeuchi \cite{Tak:GAAA}, following
work of Sweedler \cite{Swe:GSA}
in which $R$ is by assumption commutative. Actually Takeuchi invents
his $\tR$-bialgebras from different motivations, involving generalizations
of Brauer groups, and does not seem to be thinking of groupoids at all.
Lu \cite{Lu:HAQG} and Xu \cite{Xu:QG} reinvent his notion, now with
the motivation by noncommutative-geometric groupoids in mind.
(Actually most of Lu's or Xu's definition is the very same as 
Takeuchi's up to changes in notation, at least as far as 
comultiplication is concerned. For a detailed translation, and the
removal of any doubt about the notion of counit, consult
Brzezi\'nski and Militaru \cite{BrzMil:BBD}.) 

Of the other possible definitions of a quantum groupoid we should
mention the weak Hopf algebras of B\"ohm and Szlach\'anyi
\cite{BohSzl:CCQGND}, see also the recent survey 
\cite{NikVai:FQGA} by Nikshych and Vainerman 
and the literature cited there, and the notion
of a face algebra due to Hayashi \cite{Hay:FAIGQGT,Hay:BIFA}.
Face algebras were shown to be precisely the $\tR$-bialgebras
in which $R$ is commutative and separable in \cite{Sch:FATRB}.
Also, face algebras are precisely the weak bialgebras whose
target counital subalgebras are commutative. Etingof and Nikshych
\cite{EtiNik:DQGR1}
have shown that weak Hopf algebras are $\tR$-bialgebras. In fact
weak bialgebras are precisely those $\tR$-bialgebras in which 
$R$ is Frobenius-separable (for example semisimple over the
complex numbers) \cite{Sch:WHAQG}.

In the present paper we will discuss a construction that allows us
to replace the algebra $R$ in any $\tR$-bialgebra $L$
by a Morita-equivalent
algebra $S$ to obtain a $\tS$-bialgebra that has the same representation
theory, more precisely a monoidal category of representations 
equivalent to that of $L$. In fact we can, more generally, replace
$R$ by any \sqm\ equivalent algebra $S$. The notion of \sqm\ equivalence
is due to Takeuchi \cite{Tak:MT}. Two algebras $R$, $S$ are by
definition \sqm\ equivalent if we have an equivalence of $k$-linear
monoidal categories $\BiMod R\cong\BiMod S$. The definition is already
at the heart of our application: A $\tR$-bialgebra can be characterized
as having a monoidal category of representations with tensor 
product based on the tensor product in $\BiMod R$.
However, for some purposes it does seem that Morita base change
(replacing $R$ by a Morita equivalent algebra $S$) is more well-behaved
than the more general
\sqm\ base change (replacing $R$ by a \sqm\ equivalent algebra
$S$): We will show that Morita base change respects duality.

Morita (or \sqm) base change can serve two immediate purposes: One
is to produce new examples of quantum groupoids. The other, and
perhaps more useful one, is to on the contrary reduce the supply
of essentially different examples --- we can consider two 
$\tR$-bialgebras to be not very essentially different if they
are obtained from each other by \sqm\ base change. Note that the
equivalence
relation thus imposed on $\tR$-bialgebras is weaker than the
natural relation that would consider two $\tR$-bialgebras to be
equivalent if their monoidal categories of representations are
equivalent. In fact this latter equivalence relation is known to
be important and nontrivial also in the realm of ordinary bialgebras,
where Morita base change is meaningless. Thus Morita base change
presents a possibility of relating different 
$\tR$-bialgebras very closely, in a way that cannot occur between
ordinary bialgebras.

Let us state very briefly two ways in which \sqm\ base change 
reduces the supply of examples: If $R$ is an Azumaya $k$-algebra,
then any $\tR$-bialgebra is, up to \sqm\ base change, an ordinary
bialgebra. Over the field of complex numbers, every weak bialgebra
is, up to Morita base change, a face algebra.
Of course, in neither case our results show that certain 
$\tR$-bialgebras are entirely superfluous, since examples may 
occur in natural situations that come with a specific choice of
$R$.

The plan of the paper is as follows:
After recalling some definitions in \nmref{sec:tRb} and
\nmref{sec:sqm} we present the general \sqm\ base change procedure
in \nmref{sec:sqmbc}. More detailed information on 
Morita base change will be given in \nmref{sec:mbc}.
In \nmref{sec:ctd} we discuss the canonical Tannaka duality
of Hayashi \cite{Hay:CTDFSTC,Hay:BIFA}; this construction assigns
a face algebra $F$ to any finite split semisimple $k$-linear 
monoidal category. For example, it assigns such a face algebra
to the category of representations of a split semisimple 
(quasi)Hopf algebra $H$. At first sight, there is no apparent relation
between the original $H$ and Hayashi's $F$ (beyond, of course, the
fact that their monoidal representation categories are equivalent).
We show that $F$ can be obtained from $H$ in two steps: First, one
applies a kind of smash product construction that builds from 
$H$ a $\times_H$-bialgebra isomorphic to $H\o H\o H^*$ as a vector
space. Next, applying Morita base change to replace the base $H$
by the Morita equivalent product of copies of the field, one obtains
a face algebra --- which turns out to be
Hayashi's face algebra $F$.
In \nmref{sec:II1} we compute the dimension of the face algebra
obtained by Morita base change from a certain weak Hopf algebra 
constructed by Nikshych and Vainerman from a subfactor of a 
type $\II_1$ factor. It turns out that Morita base change reduces the
dimension from $122$ to $24$ without affecting the monoidal category
of representations.

{\bf Acknowledgements:} The author is indebted to Leonid Vainerman
for interesting discussions, and in particular for some help in 
understanding the example underlying \nmref{sec:II1}.

\section{Hopf algebroids}\nmlabel{Section}{sec:tRb}

In this section we will very briefly recall 
the necessary definitions and notations
on $\tR$-bialgebras. For more details we refer to 
\cite{Swe:GSA,Tak:GAAA,Sch:BNRSTHB}. 

Throughout the paper, $k$ denotes a commutative base ring, and
all modules,  algebras,  unadorned tensor products etc.\ are 
understood to be over $k$

Let $R$ be a $k$-algebra. 
We denote the opposite algebra by $\ol R$, we let
 $R\ni r\mapsto \overline r\in \overline R$ denote the obvious 
$k$-algebra antiisomorphism, and abbreviate the enveloping
algebra $R^e:=R\o\ol R$. We write $r\ol s:=r\o\ol s\in R\o\ol R$
for $r,s\in R$.

For our purposes, a handy characterization of $\tR$-bialgebras is 
the following \cite[Thm.5.1]{Sch:BNRSTHB}: A $\tR$-bialgebra $L$ is
an $R^e$-ring (that is, a $k$-algebra equipped with a 
$k$-algebra map $R^e\rightarrow L$, which we write
$r\o\ol s\mapsto r\ol s$) for which the category $\LMod L$
is equipped with a monoidal category structure such that the
``underlying'' functor $\LMod L\rightarrow \LMod{R^e}$ is a 
strict monoidal functor. Here, the monoidal category structure
on $\LMod{R^e}$ is induced via the identification with the
category $\BiMod R$ of bimodules; we denote tensor product
in $\LMod{R^e}$ by $\nti R$, or $\nt$ if no confusion is likely.

Thus, for two $L$-modules $M,N$, there is an $L$-module structure
on $M\nti RN$, and this tensor product of $L$-modules defines
a monoidal category structure on $\LMod L$.
The connection with the original definition in  \cite{Tak:GAAA}
is that the module structure on $M\nt N$ can be described in terms
of a certain comultiplication on $L$, which, however, has a more
intricate definition than in the ordinary bialgebra case. 
First of all, the comultiplication is an algebra map 
$L\rightarrow L\tR L$ into a certain subset $L\tR L\subset L\nti RL$
which has an algebra structure induced by that of $L\o L$, and
whose definition we shall now recall.

The notations $\int_r:=\int_{r\in R}$ and $\int^r:=\int^{r\in R}$,
which we will introduce only by example, are due to
MacLane, see \cite{Swe:GSA,Tak:GAAA}.
For $M,N\in\Bimod{R^e}{R^e}$ we let
$$  \int_r {_{\ol r}M}\o {_rN}:=M\o N\big/\langle \ol rm\o n-m\o rn|r\in R,m\in M,n\in N\rangle$$
and we let
$\displaystyle  \int^r M_{\ol r}\o N_r\subset M\o N$
denote the $k$-submodule consisting of all elements
$\lsum m_i\o n_i\in M\o N$ satisfying 
$\lsum m_i\ol r\o n_i=\lsum m_i\o n_ir$ for all 
$r\in R$. 
Note $\int_r {_{\ol r}M}\o {_rN}=M\nti RN$ for $M,N\in\LMod{R^e}$.

For two $R^e$-bimodules $M$ and $N$ we abbreviate
$$M\tR N:=\int^s\int_r {_{\ol r}M_{\ol s}}\o {_rN_s}.$$ 
If $M,N$ are $R^e$-rings, then so is $M\tR N$, 
with multiplication given by
$(\sum m_i\o n_i)(\sum m'_j\o n'_j)=\sum m_im'_j\o n_in'_j$,
and $R^e$-ring structure
$$R^e\ni r\o\ol s\mapsto r\o\ol s\in M\tR N.$$

For $M,N,P\in\Bimod{R^e}{R^e}$ one defines
$$M\tR P\tR N:=\int^{s,u}\int_{r,t}{_{\ol r}M_{\ol s}}
   \o{_{r,\ol t}P_{s,\ol u}}\o{_tN_u}$$
(where $\int^{s,u}:=\int^s\int^u=\int^u\int^s$). There 
are associativity maps
\begin{align*}
  (M\tR P)\tR N &\overset\alpha\rightarrow M\tR P\tR N\\
  M\tR (P\tR N) &\overset{\alpha'}\rightarrow M\tR P\tR N
\end{align*}
given on elements by the obvious formulas (doing nothing), but which 
need not be isomorphisms. If $M,N$ and $P$ are $R^e$-rings, so
is $M\tR N\tR P$, and $\alpha,\alpha'$ are $R^e$-ring maps.

An $R^e$-ring structure on the algebra
$\E=\End(R)$ is given by 
$r\o\ol s\mapsto(t\mapsto rts)$.
We have, for any $M\in\Bimod{R^e}{R^e}$, two $R^e$-bimodule
maps
\begin{xalignat*}2
  \theta :M\tR\End(R)&\rightarrow M;& m\o f&\mapsto \ol{f(1)}m\\
  \theta' :\End(R)\tR M&\rightarrow M;& f\o m&\mapsto f(1)m.
\end{xalignat*}
which are $R^e$-ring homomorphisms if $M$ is an $R^e$-ring.

Now we are prepared to write down the definition of a
$\tR$-bialgebra $L$. This is by definition an $R^e$-ring
equipped with a comultiplication, a 
map $\Delta\colon L\rightarrow L\tR L$ of $R^e$-rings
over $R^e$, and a counit, a map
$\epsilon:L\rightarrow \E$ of $R^e$-rings,
such that 
\begin{gather}
  \label{tRcoass}\alpha(\Delta\tR L)\Delta=\alpha'(L\tR\Delta)\Delta\colon L\rightarrow L\tR L\tR L\\
  \label{tRcounit}
  \theta(L\tR\epsilon)\Delta=\id_L=\theta(\epsilon\tR L)\Delta.
\end{gather}

For $\tR$-bialgebras we will make use of the usual
Sweedler notation, 
writing $\Delta(\ell)=:\ell\sw 1\o\ell\sw 2\in L\tR L$.

If $L$ is a $\tR$-bialgebra, then 
the module structure on the tensor product $M\nti RN$
of $M,N\in\LMod L$ can be described in terms of the comultiplication
of $L$ by the usual formula
$\ell (m\o n)=\ell\sw 1 m\o\ell\sw 2 n$.

The suitable definition of comodules over a $\tR$-bialgebra $L$ 
is as follows:
A left $L$-comodule is an $R$-bimodule
  $M$ together with a map 
  $\lambda\colon M\rightarrow L\tR M$ of $R$-bimodules
  such that
$$\alpha'(L\tR\lambda)\lambda=\alpha(\Delta\tR M)\lambda\colon M\rightarrow L\tR L\tR M$$
and $\theta'(\epsilon\tR M)\lambda=\id_M$ hold.
We will denote by $\LComod L$ the category of left 
$L$-comodules. We will use Sweedler notation in the form
$\lambda(m)=m\sw{-1}\o m\sw 0$ and
$\alpha(\Delta\tR M)(m)=m\swm 2\o m\swm 1\o m\sw 0$
for $L$-comodules.

The category $\LComod L$ of left $L$-comodules over a $\tR$-bialgebra
is monoidal. The tensor product of $M,N\in\LComod L$ is their
tensor product $M\ou RN$ over $R$, equipped with the comodule
structure
  \begin{align*}
    M\ou RN&\rightarrow L\tR(M\ou RN)\\
    m\o n&\mapsto m\sw{-1}n\sw{-1}\o m\sw 0\o n\sw 0
  \end{align*}

In \cite[Thm.\ and Def.3.5]{Sch:DDQG} we have introduced a notion
of $\tR$-Hopf algebra. It is rather different from that of a 
Hopf algebroid given by Lu \cite{Lu:HAQG}, although Lu's
bialgebroids are the same as $\tR$-bialgebras.
By definition, a $\tR$-bialgebra is a $\tR$-Hopf algebra
if and only if the map
$$\beta\colon L\ou{\ol R}L\ni \ell\o m\mapsto \ell\sw 1\o\ell\sw 2m\in L\nt L$$
is a bijection. Note that this is equivalent to a well-known 
characterization of Hopf algebras among ordinary bialgebra.
By \cite{Sch:DDQG} it is equivalent to saying that the underlying
functor $\LMod L\rightarrow\LMod{R^e}$ preserves inner hom-functors. 
More precisely, for each $M\in\LMod L$ the functor 
$\LMod L\ni N\mapsto N\nt M\in\LMod L$ has a right adjoint
$\hom(M,\leer)$. The $\tR$-bialgebra is a $\tR$-Hopf algebra
if and only if a canonically defined map
$\hom(M,N)\rightarrow\Hom_{\ol R-}(M,N)$ is a bijection for
all $M,N\in\LMod L$.

Let us finally recall two special cases of the notion of a 
$\tR$-bialgebra. Weak bialgebras and weak Hopf algebras
were introduced by B\"ohm and Szlach\'anyi 
\cite{BohSzl:CCQGND}. We refer to the survey \cite{NikVai:FQGA}
by Nikshych and Vainerman and the literature cited there.
It was shown in \cite{EtiNik:DQGR1} that weak Hopf algebras
are $\tR$-bialgebras. More details and a converse are in 
\cite{Sch:WHAQG}.
By definition, a weak bialgebra $H$ is a $k$-coalgebra and 
$k$-algebra
such that comultiplication is multiplicative, but not necessarily
unit-preserving (and neither is multiplication assumed to be
comultiplicative). There are specific axioms replacing the
``missing'' compatibility axioms for a bialgebra, namely,
for $f,g,h\in H$:
\begin{gather*}
  \epsilon(fgh)=\epsilon(fg\sw 1)\epsilon(g\sw 2h)=\epsilon(fg\sw 2)\epsilon(g\sw 1h),\\
  1\sw 1\o 1\sw 2\o 1\sw 3=1\sw 1\o 1\sw 21'\sw 1\o 1'\sw 2=1\sw 1\o 1'\sw 11\sw 2\o 1'\sw 2.
\end{gather*}
If $H$ is a weak bialgebra, then 
the target counital subalgebra $H_t$ consists by definition of
all elements of the form $\epsilon(1\sw 1h)1\sw 2$ with $h\in H$.
The source counital subalgebra $H_s$ is the target counital subalgebra
in the coopposite of $H$.
It turns out that $H_t$ is a subalgebra which is Frobenius-separable
(i.e.\ a multi-matrix algebra when $k=\mathbb C$ is the field of complex
 numbers),
anti-isomorphic to $H_s$,
and that $H$ has the structure of a $\tR$-bialgebra for $R=H_t$,
in which $H_s$ is the image of $\ol R$ in $H$. Moreover, any
$\tR$-bialgebra in which $R$ is Frobenius-separable can be obtained
in this way from a weak bialgebra. A weak Hopf algebra is by definition
a weak bialgebra $H$
with an antipode, which in turn is an anti-automorphism
of $H$ whose axioms we shall not recall. The antipode maps
$H_t$ isomorphically onto $H_s$ and vice versa. We have shown
in \cite{Sch:WHAQG} that a weak bialgebra has an antipode if and
only if the associated $\tR$-bialgebra is a $\tR$-Hopf algebra.

The face algebras introduced earlier by Hayashi 
\cite{Hay:FAIGQGT,Hay:BIFA} are recovered as a yet more special case
of weak bialgebras, namely that where the target (and source)
counital subalgebra is commutative. 
In particular, as shown in \cite{Sch:FATRB}, a face algebra $H$ 
the same thing as a $\tR$-bialgebra in which $R$ is commutative
and separable.
We will only be using the case where the base
field is the field of complex numbers, so that $R$ is a direct product
of copies of the field. In particular, the images of the minimal
idempotents of $R$ in $H$ form a distinguished family of idempotents
in $H$, which feature prominently in Hayashi's original definition
(along, of course, with the images of the corresponding idempotents
in $\ol R$). We shall refer to them as the face idempotents
of $H$; their number, or the dimension of $R$, is an important 
structure element of $H$. 

\section{Morita-\ and \sqm-equivalence}\nmlabel{Section}{sec:sqm}
Suppose $R,S$ are Morita equivalent $k$-algebras (we shall write
$R\msim S$ for short). Then by definition the categories
$\LMod R$ and $\LMod S$ of left modules are equivalent as 
$k$-linear abelian categories. In this situation, one also gets
an equivalence $(\BiMod R,\ou R)\cong(\BiMod S,\ou S)$ of 
$k$-linear monoidal categories. This can be seen by applying
Watts' theorem \cite{Wat:ICSAF}, which says that the monoidal
category $\BiMod R$ can be viewed as the category of right
exact $k$-linear endofunctors of $\LMod R$. Somewhat more
useful is the following explicit description of the monoidal
equivalence: When $R\msim S$, fix a strict Morita context
$(R,S,P,Q,f,g)$. In particular, we have
$P\in\Bimod SR$, $Q\in\Bimod RS$,
$f\colon P\ou RQ\rightarrow S$ an isomorphism of 
$S$-bimodules and $g\colon Q\ou SP\rightarrow R$ 
an isomorphism of $R$-bimodules. 
An equivalence is given by 
$$\mathcal F\colon\LMod R\ni M\mapsto P\ou RM\in\LMod S,$$ 
and we can describe a matching equivalence of bimodule
categories
$$(\hat{\mathcal F},\xi)\colon(\BiMod R,\ou R)\rightarrow(\BiMod S,\ou S)$$
as follows:
We set
$\hat{\mathcal F}(M)=P\ou RM\ou RQ$ as $S$-bimodules, and we define
the monoidal functor structure 
$$\xi\colon\mathcal F(M)\ou S\mathcal F(N)\rightarrow\mathcal F(M\ou RN)$$
as the composition
\begin{multline*}
P\ou RM\ou RQ\ou SP\ou RN\ou RQ\xrightarrow{P\o M\o g\o N\o Q}
P\ou RM\ou RR\ou RN\ou RQ\\\cong P\ou RM\ou RN\ou RQ
\end{multline*}
It is useful to know that the equivalence $\mathcal F$ and the
monoidal equivalence $\hat{\mathcal F}$ are compatible in the 
following sense: The category $\LMod R$ is in a natural way
a left $\BiMod R$-category in the sense of Pareigis \cite{Par:NARMTII},
that is, a category on which $\BiMod R$ acts (by tensor product).
The compatibility says that the following diagram commutes
up to coherent natural isomorphisms:
$$\xymatrix{\BiMod R\times\LMod R\ar[r]^-{\ou R}\ar[d]^{\hat\F\times\F}
&\LMod R\ar[d]^{\F}\\
\BiMod S\times\LMod S\ar[r]^-{\ou S}&\LMod S}$$

Takeuchi \cite{Tak:MT} has introduced and investigated the 
notion of \sqm-equivalence of $k$-algebras; by his definition,
two $k$-algebras $R,S$ are \sqm-equivalent, written
$R\sqmsim S$, if there is an equivalence of $k$-linear monoidal 
categories $\BiMod R\cong\BiMod S$. By the above, 
Morita equivalence clearly implies \sqm-equivalence. On the other
hand, since $\BiMod R\cong\LMod{R^e}$, the enveloping algebras
of \sqm-equvialent algebras are Morita equivalent, so that
$$R\msim S\quad\Rightarrow\quad R\sqmsim S\quad\Rightarrow\quad R^e\msim S^e.$$
Neither of the reverse implications holds.

Note that a bimodule $M\in\BiMod R$ has a left dual object in the
monoidal category $(\BiMod R,\ou{R})$ if and only if it is 
finitely generated projective as a right $R$-module. It follows that
any equivalence of monoidal categories $\BiMod R\cong\BiMod S$
maps left (or right) finitely generated projective modules to
left (or right) finitely generated projective modules. This reduces
to a standard fact on projective modules, if the equivalence comes
from a Morita equivalence $R\msim S$, for then it maps
$M\in\BiMod R$ to $P\ou RM\ou RQ$, which is finitely generated 
projective as left $S$-module if $_RM$ is finitely generated projective,
for $_SP$ and $_RQ$ are finitely generated projective.

\section{A \sqm-base change principle}\nmlabel{Section}{sec:sqmbc}
An $A$-ring $L$ for a $k$-algebra $A$ is an algebra in the monoidal
category of $A$-bimodules. As an immediate consequence, if
$A\msim B$, then an $A$-ring is essentially the same as a 
$B$-ring. Moreover, if $\wtilde L$ is the $B$-ring corresponding to
the $A$-ring $L$, then $\wtilde L$-modules are essentially the 
same as $L$-modules, since the actions of $\BiMod B$ on $\LMod B$
and of $\BiMod A$ on $\LMod A$ are compatible with the equivalences.
Thus we have:
\begin{Lem}
  Let $L$ be an $A$-ring over the $k$-algebra $A$.
  Assume given a $k$-algebra $B$ and a strict Morita context
  $(A,B,C,D,\phi,\psi)$. Then there is a $B$-ring $\wtilde L$, and
  a category equivalence 
  $\mathcal G\colon\LMod{L} \rightarrow\LMod{\wtilde L}$ 
  lifting the equivalence
  $\F\colon\LMod A\rightarrow \LMod B$ given by tensoring with $C$, as 
  in the following diagram:
  $$\xymatrix{\LMod L\ar[r]^{\mathcal G}\ar[d]
     &\LMod{\wtilde L}\ar[d]\\\LMod A\ar[r]^{\mathcal F}&\LMod B}$$
  in which the vertical arrows are the underlying functors induced
  by the $A$-ring structure of $L$, and the $B$-ring structure
  of $\wtilde L$, respectively.
\end{Lem}
Explicitly, the $B$-ring $\wtilde L $ is given by
$\wtilde L :=C\ou AL \ou AD$ with unit map
$$B\cong C\ou AD\xrightarrow{C\ou A\eta\ou AD}C\ou AL \ou AD$$
in which $\eta$ is the unit map of the $A$-ring $L $, and 
multiplication map
\begin{multline*}
\wtilde L \ou B\wtilde L =C\ou AL \ou AD\ou BC\ou AL \ou AD
\cong C\ou AL \ou AL \ou AD\xrightarrow{C\ou A\nabla\ou AD}
C\ou AL \ou AD.\end{multline*}
When $M$ is a $L $-module, then the $\wtilde L $-module $\F(M)$
is $C\ou AM$ equipped with the $\wtilde L $-module structure
$$\wtilde L \ou B\F(M)=C\ou AL \ou AD\ou BC\ou AM=C\ou AL \ou AM
\xrightarrow{C\ou A\mu_M}C\ou AM=\F(M).$$
\begin{Rem}\nmlabel{Remark}{Brattelrem}
  Assume that the algebras involved in the situation 
  above are multi-matrix algebras over a
  field $k$.
 
  Then the Bratteli diagram of the inclusion $B\subset \wtilde L$
  is the same as that of the inclusion $A\subset L$, except that
  the ranks of the components of $A$ have to be replaced by
  the ranks of the components of $B$, and the top floor 
  representing $\wtilde L$ has to be adjusted accordingly.
\end{Rem}
In fact the number of edges between a vertex on the bottom floor
of the Bratelli diagram and a vertex on the top floor is the 
multiplicity of an irreducible representation of $A$ in a certain
irreducible representation of $L$. Evidently these numbers do not
change in the construction in the Lemma.
  
\begin{Thm}\nmlabel{Theorem}{sqmbcthm}
  Let $L$ be a $\times_R$-bialgebra for a $k$-algebra $R$. 
  Let $S$ be a $k$-algebra which is \sqm-equivalent to $R$. Then 
  there is a $\times_S$-bialgebra $\wtilde L$ whose module category
  is equivalent to that of $L$, as a monoidal category.

  More precisely, assume given a monoidal category equivalence
  $\F\colon \LMod{R^e}\rightarrow\LMod{S^e}$. Then there is
  a $\times_S$-bialgebra $\wtilde L$ and a monoidal category 
  equivalence $\mathcal G\colon\LMod L\rightarrow\LMod{\wtilde L}$
  making 
  \begin{equation}\label{monfundiag}
    \begin{array}{c}
    \xymatrix{\LMod L\ar[r]^{\mathcal G}\ar[d]^{\mathcal U}
     &\LMod{\wtilde L}\ar[d]^{\wtilde{\mathcal U}}
     \\\LMod{R^e}\ar[r]^{\F}&\LMod{S^e}}
    \end{array}
  \end{equation}
  a commutative diagram of monoidal functors (in which the vertical
  arrows are underlying functors.
\end{Thm}
\begin{proof}
  We know already that there is an $S^e$-ring $\wtilde L$ and a 
  category equivalence $\mathcal G$ making the
  square in the Theorem a commutative diagram of $k$-linear functors.
  We can endow $\LMod{\wtilde L}$ with a monoidal category structure 
  such that $\mathcal G$ is a monoidal functor. 
  Then the underlying functor
  $\wtilde{\mathcal U}$ is monoidal as well, since it can be written as
  the composition $\wtilde{\mathcal U}=\F\mathcal U\mathcal G\inv$.
  Now \cite[Thm.5.1]{Sch:BNRSTHB} implies that there exists a 
  $\times_S$-bialgebra structure on $\wtilde L$ inducing the given
  monoidal category structure on $\LMod{\wtilde L}$.
\end{proof}
\begin{Def}
  Let $L $ be a $\times_R$-bialgebra, and $S$ a 
  $k$-algebra \sqm-equivalent to $R$. We will say that the 
  $\times_S$-bialgebra $\wtilde L$ obtained from $L $ as in the
  proof of \nmref{sqmbcthm} is obtained from $L $ by a 
  \sqm\ base change.
\end{Def}
Thus a $\times_S$-bialgebra obtained from a $\times_R$-bialgebra
$L $ by a \sqm-base change has the same monoidal representation
category as $L $ itself. The difference is ``merely'' in a change
of the base algebra.

We will be somewhat sloppy in our terminology: Given a $\tR$-bialgebra
$L$ and $S\sqmsim R$, we will speak of {\em the\/} $\tS$-bialgebra
$\wtilde L$ obtained from $L$ by \sqm\ base change. This suppresses
the choice of a monoidal category equivalence $\BiMod R\cong\BiMod S$,
by which (and not by $S$ alone) $\wtilde L$ is determined.

\begin{Cor}\nmlabel{Corollary}{Azuchange}
  Let $L $ be a $\times_R$-bialgebra, where $R$ is an 
  Azumaya $k$-algebra. Then $L $ can be obtained by \sqm-base
  change from an ordinary bialgebra $H$.
\end{Cor}
\begin{proof}
  Takeuchi \cite[Ex.2.4]{Tak:MT} has shown that the algebra $R$ is
  Azumaya if and only if $R\sqmsim k$.
\end{proof}
Takeuchi also gives the following elegant description of the 
monoidal category equivalence $\BiMod R\rightarrow\kmod$:
It maps $M$ to the centralizer $M^R$ of $R$ in $M$, whereas its
inverse maps $V\in\kmod$ to $V\o R$ with the obvious $R$-bimodule
structure. Thus, the ordinary $k$-bialgebra associated to a 
$\times_R$-bialgebra $L $ is
$$H=\{\ell \in L |\forall r\in R\colon r\ell=\ol r\ell\wedge\ell r=\ell\ol r\}$$
whereas $L$ can be obtained from $H$ by merely tensoring with
two copies of $R$, one of which gives the left $R^e$-module structure, 
and the other one the right $R^e$-module structure of
$R\o L\o R$. 

In a sense our corollary says that examples of $\tR$-bialgebras
in which $R$ is Azumaya are irrelevant; they are just versions
of ordinary bialgebras in which the base ring is enlarged, without
affecting the representation theory.
That would also apply to examples like that considered by Kadison
in \cite[Thm.5.2]{Kad:HAHSE}. In fact in the example of a 
$\times_R$-bialgebra $T$ there, $R$ is Azumaya over its center 
$Z$. Moreover, the two algebra maps 
$Z\rightarrow T$ coming from the maps $R\rightarrow T$ and
$\ol R\rightarrow H$ coincide by construction and have central
image. Thus $T$ can be considered as a $\times_R$-bialgebra for
the $Z$-algebra $R$; since this is Azumaya, \nmref{Azuchange}
applies, so that $T$ can be obtained by \sqm\ base change from 
a $Z$-bialgebra. However, we should rush to concede that 
\nmref{Azuchange} does of course not rule out that 
interesting examples of $\tR$-bialgebras over Azumaya
$k$-algebras arise naturally. In fact Kadison's 
example is constructed from a natural situation
that comes with a natural choice of base $R$. 
Moreover, the example gives us
the opportunity to point out a 
certain subtlety about \sqm\ base change: The $\tR$-bialgebra
$T$ in \cite[Sec.4]{Kad:HAHSE} occurs in duality with another
$\tR$-bialgebra $S$, which can also be considered as a $\tR$-bialgebra
for the $Z$-algebra $R$. Thus, if $R$ is Azumaya over $Z$, then 
$S$ can be reduced by \sqm\ base change to a $Z$-bialgebra $S'$,
while $T$ can be replaced by a $Z$-bialgebra $T'$. However, we
do not have any indication that $S'$ and $T'$ are still dual
to each other. It is conceivable that the duality only shows over
the ring $R$. We will show below that Morita base change is 
compatible with duality.

Closing the section, let us show that \sqm\ base change preserves
the property of a $\tR$-bialgebra of being a $\tR$-Hopf algebra
in the sense of \cite{Sch:DDQG}:
\begin{Prop}
  Let $L$ be a $\tR$-bialgebra, $S\sqmsim R$, and let $\wtilde L$
  be the $\tS$-bialgebra obtained from $L$ by \sqm\ base change.

  $\wtilde L$ is a $\tS$-Hopf algebra if and only if $L$ is a 
  $\tR$-Hopf algebra.
\end{Prop}
\begin{proof}
  In the diagram \eqref{monfundiag}, the horizontal functors are
  monoidal equivalences, hence preserve inner hom-functors.
  Thus the left hand vertical functor preserves inner hom-functors
  if and only the right hand one does. 
\end{proof}

\section{Morita base change}\nmlabel{Section}{sec:mbc}
Morita equivalence implies \sqm\ equivalence. Thus, given
a $\tR$-bialgebra $L$ and a $k$-algebra $S$ Morita equivalent
to $R$, we can apply \sqm\ base change (which, of course, we
shall call Morita base change in this case) to $L$ to obtain
a $\times_S$-bialgebra $\wtilde L$ with equivalent monoidal 
module category.
\subsection{Morita base change --- explicitly}
To find out what the result looks like more explicitly, fix a 
Morita context $(R,S,P,Q,f,g)$.
We will write
$f(p\o q)=pq$, $g(q\o p)=qp$,
$f\inv(1_S)=p_i\o q^i\in P\ou R Q$ and $g\inv(1_R)=q_i\o p^i\in Q\ou S P$
(with a summation over upper and lower indices understood).
Write $\ol P\in\Bimod{\ol R}{\ol S}$ for the bimodule opposite
to $P$, and $\ol p$ with $p\in P$ for a typical element; similarly
for $\ol Q\in\Bimod{\ol S}{\ol R}$. Somewhat dangerously we write
$P^e:=P\o\ol Q\in\Bimod{S^e}{R^e}$ and $Q^e:=Q\o\ol P\in\Bimod{R^e}{S^e}$,
so that the bimodules $P^e$ and $Q^e$ induce the equivalence
$\LMod{R^e}\cong\LMod{S^e}$ underlying the \sqm\ equivalence between
$R$ and $S$ induced by the Morita equivalence between $R$ and $S$.

To keep our formulas a manageable size, we will write
$p\ol q:=p\o\ol q\in P\o\ol Q=P^e$, and similarly for the typical
elements of $Q^e$.

Now let $L$ be a $\tR$-bialgebra. The $\tS$-bialgebra $\wtilde L$
obtained from 
$L$ by Morita base change has underlying $S^e$-bimodule
$P^e\ou{R^e}L\ou{R^e}Q^e$. 
The equivalence $\LMod L\cong\LMod{\wtilde L}$ sends 
$M\in\LMod L$ to $P^e\ou{R^e}M$, with the $\wtilde L$-module
structure given by
$$(p_1\ol{q_1}\o\ell\o q_2\ol{p_2})(p_3\ol{q_3}\o m)
  =p_1\ol{q_1}\o\ell(q_2p_3)(\ol{q_3p_2})m$$
for $p_1,p_2,p_3\in P$ and $q_1,q_2,q_3\in Q$.
The monoidal functor structure of the equivalence is given by
\begin{align*}
\xi\colon(P^e\ou{R^e}M)\nti S(P^e\ou{R^e}N)&\wtilde{\rightarrow} P^e\ou{R^e}(M\nti R N)\\
 p_1\ol{q_1}\o m\o p_2\ol{q_2}\o n&\mapsto p_1\ol{q_2}\o m\o(q_1p_2)n\\
                                      &=p_1\ol{q\sw 2}\o\ol{q_1p_2}m\o n\\
p\ol{q_i}\o m\o p^i\ol q&\mapsfrom p\ol q\o m\o n
\end{align*}
for $M,N\in\LMod L$, $m\in M$, $n\in N$, $p_1,p_2\in P$, $q_1,q_2\in Q$.

It follows that the $\wtilde L$-module structure of the tensor product
of two $\wtilde L$-modules coming via the equivalence from $L$-modules
$M,N$ can be computed as the composition
\begin{multline*}
  (P^e\ou{R^e}L\ou{R^e}Q^e)\ou{S^e}((P^e\ou{R^e}M)\nti S(P^e\ou R^eN))
    \\\xrightarrow{\id\o\xi}
    (P^e\ou{R^e}L\ou{R^e}Q^e)\ou{S^e}(P^e\ou{R^e}(M\nti RN))
    \\\xrightarrow{\,\mu\,}
    P^e\ou{R^e}(M\nti RN)
    \xrightarrow{\xi\inv}
    (P^e\ou{R^e}M)\nti S(P^e\ou{R^e}N)
\end{multline*}
hence is given by
\begin{multline*}
  (p_1\ol{q_1}\o\ell\o q_2\ol{p_2})(p_3\ol{q_3}\o m\o p_4\ol{q_4}\o n)
    \\=\xi\inv((p_1\ol{q_1}\o\ell\o q_2\ol{p_2})(p_3\ol{q_4}\o m\o(q_3p_4)n)
    \\=\xi\inv(p_1\ol{q_1}\o\ell(q_2p_3)(\ol{q_4p_2})(m\o(q_3p_4)n))
    \\=\xi\inv(p_1\ol{q_1}\o\ell\sw 1(q_2p_3)m\o\ell_2(\ol{q_4p_2})(q_3p_4)n)
    \\=p_1\ol{q_i}\o\ell\sw 1(q_2p_3)m\o p^i\ol{q_1}\o\ell_2(\ol{q_4p_2})(q_3p_4)n
    \\=p_1\ol{q_i}\o\ell\sw 1(q_2p_3)m\o p^i\ol{q_1}\o\ell_2(q_3p_j)(q^jp_4)(\ol{q_4p_2})n
    \\=p_1\ol{q_i}\o\ell\sw 1(q_2p_3)(\ol{q_3p_j})m\o p^i\ol{q_1}\o\ell\sw 2(q^jp_4)(\ol{q_4p_2})n
    \\=(p_1\ol{q_i}\o\ell\sw 1\o q_2\ol{p_j})(p_3\ol{q_3}\o m)
         \o(p^i\ol{q_1}\o\ell\sw 2\o q^j\ol{p_2})(p_4\ol{q_4}\o n)
\end{multline*}
for $p_1,\dots,p_4\in P$, $q_1,\dots,q_4\in Q$, $\ell\in L$, 
$m\in M$, and $n\in N$. this proves that the comultiplication in 
$\wtilde L$ is given by the formula
$$\Delta(p_1\ol{q_1}\o\ell\o q_2\ol{p_2})
    =(p_1\ol{q_i}\o\ell\sw 1\o q_2\ol{p_j})\o(p^i\ol{q_1}\o\ell\sw 2\o q^j\ol{p_2})
$$
for $p_1,p_2\in P$ and $q_1,q_2\in Q$.
\subsection{Weak bialgebras versus face algebras}
Let us be yet more concrete for the case that $R$ is a multi-matrix
algebra $R=\bigoplus_{\alpha=1}^n M_{d_\alpha}(k)$, and $S=k^n$. A Morita
context $(R,S,P,Q,f,g)$ can be given as follows:
$P$ is generated as a right $R$-module by one element $p$
which is a sum $p=\sum_{\alpha=1}^nE_{11}^{(\alpha)}$ of minimal idempotents
(where we have denoted the matrix units in the $\alpha$-th component
by $E^{(\alpha)}_{ij}).$
$Q$ is generated as left $R$-module by the same element $p$.
both maps $f,g$ are given by matrix multiplication. We have
$f\inv(1_S)=p\o p$, and
$g\inv(1_R)=\sum_{\alpha=1}^n\sum_{i=1}^{d_\alpha}E^{\alpha}_{i1}\o E^{(\alpha)}_{1i}$.
Let $L$ be a $\tR$-bialgebra, and $\wtilde L$ the $\tS$-bialgebra 
obtained from it by Morita base change. Then 
$\wtilde L=p\ol pLp\ol p\subset L$, with multiplication given by
multiplication in $L$, unit $p\ol p$, and comultiplication 
$$\Delta(p\ol p\ell p\ol p)=
    p\ol {E^{(\alpha)}_{i1}}
       \ell\sw 1p\ol{p}\o E^{(\alpha)}_{1i}\ol{p}\ell\sw 2p\ol{p}.
$$

Now let $k=\mathbb C$ be the field of complex numbers. Then 
a $\tR$-bialgebra for a multi-matrix algebra $R$ is the same as
a weak bialgebra in the sense of B\"ohm and \Szlachanyi\
\cite{BohSzl:CCQGND,BohNilSzl:WHAIITCS}. If $R$ is commutative, then this is in turn
the same thing as a face algebra in the sense of Hayashi
\cite{Hay:BIFA}. Thus Morita base change says that Hayashi's
face algebras are a sufficiently general case of weak bialgebras,
at least as long as we are interested in the respective 
module categories:

\begin{Cor}
  Let $H$ be a weak bialgebra over the field of complex numbers.

  Then $H$ can be obtained by Morita base change from a weak
  bialgebra whose source counital subalgebra is commutative.

  In particular, 
  there is a face algebra $F$ and a monoidal category
  equivalence $\LMod H\cong\LMod{F}$.

  $H$ is a weak Hopf algebra if and only if $F$ is a face Hopf algebra.
\end{Cor}

\begin{Rem}
  For the case of semisimple $H$, it follows in fact from Hayashi's
  canonical Tannaka duality \cite{Hay:BIFA,Hay:CTDFSTC} that there
  is a face algebra $F$ and a monoidal category equivalence
  $\LMod H\cong\LMod{F}$. The corollary above shows the same 
  for non-semisimple $H$, but it is also a different result in  
  the semisimple case. A peculiar feature of Hayashi's
  canonical Tannaka duality (on which we will give more details
  in the next section) is that it yields semisimple face algebras
  with the same number of face idempotents and irreducible 
  representations. This clearly 
  needs not be the case for the face algebras obtained by 
  Morita base change.
  A trivial example is the trivial Morita
  base change applied to an ordinary Hopf algebra, which leaves
  us with the same Hopf algebra, or only one face idempotent. 
  More examples will appear below.
\end{Rem}
\subsection{Duality}
There is a well-behaved notion of duality for $\tR$-bialgebras,
developed in \cite{Sch:DDQG}, and shown to be compatible with
the duality for weak bialgebras in \cite{Sch:WHAQG}. The
main difficulty in the definitions is to sort out how the four 
module structures in a $\tR$-bialgebra should be translated
through the duality, and to check that the formulas defining 
the dual structures are well defined with respect to the various
tensor products over $R$. A specialty is that one can define the
$\tR$-bialgebra analog of the opposite or coopposite of the dual
of an ordinary bialgebra, but not the direct analog of the dual
(unless one wants to allow two versions of ``left'' and ``right''
bialgebroids like Kadison and Szlach\'anyi \cite{KadSzl:DBDTRE}).
More generally, one defines 
\cite[Def.5.1]{Sch:DDQG} a skew pairing between two
$\tR$-bialgebras $\Lambda$ and $L$ to be a $k$-linear map
$\tau\colon\Lambda\o L\rightarrow R$ satisfying
  \begin{alignat}
    \tau((r\o \ol s)\xi(t\o\ol u)|\ell)v
      &=r\tau(\xi|(t\o\ol v)\ell(u\o\ol s)),
      \label{skp.1}\\
    \tau(\xi|\ell m)
       &=\tau(\ol{\tau(\xi\sw 2|m)}\xi\sw 1|\ell),
       \label{skp.2}
       &\tau(\xi|1)=\epsilon(\xi)(1),\\
    \tau(\xi\zeta|\ell)
              &=\tau(\xi|\tau(\zeta|\ell\sw 1)\ell\sw 2),
       \label{skp.3}
       &\tau(1|\ell)=\epsilon(\ell)(1)
  \end{alignat}
for all $\xi,\zeta\in\Lambda$, $\ell,m\in L$, $r,s,t,u,v\in R$.
\begin{Prop}
Let $\tau\colon\Lambda\o L\rightarrow R$ be a skew pairing between
$\tR$-bialgebras $\Lambda$ and $L$. Let $S\msim R$, and let
$\wtilde\Lambda,\wtilde L$ be the $\tS$-bialgebras obtained from 
$\Lambda$ and $L$ by Morita base change.

Then
a skew pairing $\wtilde\tau\colon\wtilde\Lambda\o\wtilde L\rightarrow S$
can be defined by
\begin{align*}
  \wtilde\tau(p_1\o\ol{q_1}\o\xi\o q_2\o\ol{p_2}| p_3\o\ol{q_3}\o\ell\o q_4\o\ol{p_4})
    &=p_1\tau(\ol{q_1p_4}\xi(q_2p_3)(\ol{q_4p_2})|\ell)q_3\\
    &=p_1\tau(\xi|(q_2p_3)\ell(q_4p_2)(\ol{q_1p_4}))q_3
\end{align*}
for $p_1,\dots,p_4\in P$, $q_1,\dots,q_4\in Q$, $\xi\in\Lambda$,
and $\ell\in L$.
\end{Prop}
\begin{proof}
To make the formulas a manageable size, we abbreviate
$p\ol q:=p\o\ol q\in P^e=P\o\ol Q$ for $p\in P$ and $q\in Q$, 
and similarly for elements of $Q^e$.

We omit checking \eqref{skp.1} for $\wtilde\tau$. 
Now let
$p_1,\dots,p_6\in P$, $q_1,\dots,q_6\in Q$,
$\xi\in\Lambda$, and $\ell,m\in L$. Then 
\begin{multline*}
\wtilde\tau((p_1\ol{q_1}\o\xi \o q_2\ol{p_2})\sw 1|
           p_3\ol{q_3}\o \ell\o  q_4\ol{p_4}
           \ol{\wtilde\tau((p_1\ol{q_1}\o\xi\o q_2\ol{p_2})\sw 2|
           p_5\ol{q_5}\o m\o q_6\ol{p_6})})
\\=
\wtilde\tau(p_1\ol{q_i}\o \xi\sw 1\o q_2\ol{p_j}|
           p_3\ol{q_3}\o\ell\o q_4\ol{p_4}
           \ol{\wtilde\tau(p^i\ol{q_1}\o\xi\sw 2\o q^j\ol{p_2}|
           p_5\ol{q_5}\o m\o q_6\ol{p_6})})
\\=
\wtilde\tau(p_1\ol{q_i}\o\xi\sw 1\o q_2\ol{p_j}|
       p_3\ol{q_3}\o \ell\o q_4\ol{p_4}
    \ol{p^i\tau(\ol{q_1p_6}\xi\sw 2(q^jp_5)(\ol{q_6p_2})|m)q_5})
\\=
\wtilde\tau(p_1\ol{q_i}\xi\sw 1q_2\ol{p_5}|p_3\ol{q_3}\ell(\ol{q_5p_4})
    \ol{\tau(\ol{q_1p_6}\xi\sw 2(\ol{q_6p_2})|m)}q_4\ol{p^i})
\\=
p_1\tau(\ol{q_ip^i}\xi\sw 1(q_2p_3)(\ol{q_4p_5})|\ell(\ol{q_5p_4})
    \ol{\tau(\ol{q_1p_6}\xi\sw 2\ol{q_6p_2}|m)})q_3
\\=
p_1\tau(\xi\sw 1(q_2p_3)|\ell(q_4p_5)
   \ol{\tau(\ol{q_1p_6}\xi\sw 2\ol{q_6p_2}|\ol{q_5p_4}m)})q_3
\\=
p_1\tau(\ol{q_1p_6}\xi(q_2p_3)(\ol{q_6p_2})|\ell(q_4p_5)(\ol{q_5p_4})m)q_3
\\=
\wtilde\tau(p_1\ol{q_1}\o\xi\o q_2\ol{p_2}|
  p_3\ol{q_3}\o\ell(q_4p_5)(\ol{q_5p_4})m\o q_6\ol{p_6})
\\=
\wtilde\tau(p_1\ol{q_1}\o\xi\o q_2\ol{p_2}|
           (p_3\ol{q_3}\o\ell\o q_3\ol{p_4})(p_5\ol{q_5}\o m\o q_6\ol{p_6}))
\end{multline*}
proves \eqref{skp.2} for $\wtilde\tau$ (we omit treating the second
part).

The proof for \eqref{skp.3} is similar.
\end{proof}
By definition, a skew pairing $\tau\colon \Lambda\o L\rightarrow R$
induces a map
$\phi\colon\Lambda\rightarrow\Hom_{\ol R-}(L,R)$.
There is an $R^e$-ring structure
\cite[Lem.5.5]{Sch:DDQG} 
on $L^\vee:=\Hom_{\ol R-}(L,R)$ 
for which $\phi$
is a morphism of $R^e$-rings. In particular, the induced
$R^e$-bimodule structure \cite[Def.5.4]{Sch:DDQG}
satisfies
$(r\ol s\xi t\ol u)(\ell)=r\xi(t\ell u\ol s)$ 
for $r,s,t,u\in R$, $\xi\in L^\vee$, and $\ell\in L$.

If $L$ is finitely generated projective as left $\ol R$-module,
then $L^\vee$ has a $\tR$-bialgebra structure 
\cite[Thm.5.12]{Sch:DDQG} such that
evaluation defined a skew pairing between $L^\vee$ and $L$.
We call this $\tR$-bialgebra the left dual of $L$.
\begin{Prop}
  Let $L$ be a $\tR$-bialgebra that is finitely generated projective
  as left $\ol R$-module. 

  Let $S\msim R$, and let $\wtilde L$ be the $\tS$-bialgebra obtained
  from $L$ by Morita base change. Then $\wtilde L$ is finitely
  generated projective as left $\ol S$-module, and its left
  dual $\tS$-bialgebra $\wtilde L^\vee$ is isomorphic to the
  $\tS$-bialgebra $\widetilde{L^\vee}$ that is obtained from the
  left dual $L^\vee$ of $L$ by Morita base change.
\end{Prop}
\begin{proof}
$\wtilde L=P^e\ou{R^e}L\ou{R^e}Q^e$ is finitely generated projective
as left $\ol S$-module since the modules $P_R$, $_{\ol S}\ol Q$,
$_{\ol R}L$, and $_{R^e}Q^e$ are finitely generated projective.

Since the $R$-modules $P$ and $Q$ are finitely generated projective
and each other's dual, we have
$\Hom(P\ou RM,V)\cong\Hom(M,V)\ou RQ$ for any $M\in\LMod R$ and
any $k$-module $V$, and similarly 
$\Hom(N\ou RQ,V)\cong P\ou R\Hom(N,V)$. 
We use this three times in the second isomorphism in the
following calculation. The first isomorphism uses the category
equivalence $\LMod{\ol R}\cong\LMod{\ol S}$ given by tensoring
with $\ol Q$. The fourth isomorphism is an instance of the
general isomorphism $\Hom_{R-}(M,P)\cong \Hom_{R-}(M,R)\ou RP$,
in the third we have used that
$\int^r$ commutes with tensor products by
flat (in particular by projective) modules.
The last step merely replaces the rightmost $P$ by $\ol P$.
\begin{multline*}
  \Hom_{\ol S-}(\wtilde L,S)
  =\Hom_{\ol S-}(P^e\ou{R^e}L\ou{R^e}Q^e,P\ou RQ)
  \cong\int^r\Hom(P\ou R{_{\ol r}L}\ou{R^e}Q^e,P_r)
  \\\cong\int^r\int_{s,t,u}P_s\o\ol Q_{\ol t}\o\Hom({_{u\ol r}L_{s\ol t}},P_r)\o{_uQ}
  \\\cong\int_{s,t,u}P_s\o\ol Q_{\ol t}\o\int^r\Hom({_{u\ol r}L_{s\ol t}},P_r)\o{_uQ}
  \\\cong\int_{s,t,u,v}P_s\o\ol Q_{\ol t}\o\Hom_{\ol R}({_uL_{s\ol t}},{_vR})\o{_uQ}\o{_vP}
  \cong P^e\ou{R^e}L^\vee\ou{R^e}Q^e.
\end{multline*}
Let $F\colon P^e\ou{R^e}L^\vee\ou{R^e}Q^e\rightarrow\wtilde L^\vee$
denote the resulting isomorphism. It is easy to check
that 
$$F(p_1\ol{q_1}\o\xi\o q_2\ol{p_2})(p_3\ol{q_3}\o\ell\o q_4\ol{p_4})
 =p_1\xi((q_2p_3)\ell(q_4p_2)(\ol{q_1p_4}))q_3$$
for $p_1,\dots,p_4\in P$, $q_1,\dots,q_4\in Q$, $\xi\in L^\vee$
and $\ell\in L$. 
In other words, the evaluation of $\wtilde L^\vee$ on $\wtilde L$
can be written as the composition
$$\wtilde L^\vee\o \wtilde L\xrightarrow{F\o\wtilde L}\widetilde{L^\vee}\o\wtilde L\xrightarrow{\wtilde\tau}S,$$
where $\wtilde\tau$ is the skew pairing of $\tS$-bialgebras
induced by the evaluation $\tau\colon L^\vee\o L\rightarrow R$.
It follows that $F$ is an isomorphism of $\tS$-bialgebras.
\end{proof}
\begin{Rem}
  Let $R,S,L$ be as above. Since $\LComod L\cong\LMod{L^\vee}$ as
  monoidal categories by \cite[Cor.5.15]{Sch:DDQG}, we can conclude
  that $\LComod L\cong\LComod{\wtilde L}$ as monoidal categories.
  Moreover, the equivalence is induced by the monoidal category
  equivalence $\BiMod R\cong\BiMod S$.
  Explicitly, it asssigns to $M\in\LComod L$ the 
  $S$-bimodule $P\ou RM\ou RQ$ endowed with the left
  $\wtilde L$-comodule structure 
  \begin{align*}
    P\ou RM\ou RQ&\rightarrow (P^e\ou{R^e}L\ou{R^e}Q^e)\tS(P\ou RM\ou RQ)\\
    p\o m\o q&\mapsto (p\o\ol{q_i}\o m\swm 1\o q\o\ol{p_j})\o(p^i\o\o m\sw 0\o q^j)
  \end{align*}
  (where $M\ni m\mapsto m\swm 1\o m\sw 0\in L\tR M$ denotes the 
  comodule structure on $M$).

  We conjecture that 
  the same formula can be used to define a category equivalence
  $\LComod L\cong\LComod{\wtilde L}$ when $L$ is not finitely generated
  projective as left $\ol R$-module. 
\end{Rem}
\begin{Rem}
  Let $L$ be a $\tR$-bialgebra which is finitely generated projective
  as left $\ol R$-module, and let $R\sqmsim S$.
  Let the monoidal equivalence $\LMod{R^e}\cong\LMod{S^e}$ be
  induced by a Morita context involving the modules
  $C\in\Bimod{S^e}{R^e}$ and $D\in\Bimod{R^e}{S^e}$. 
  By the remarks closing \nmref{sec:sqm},
  we know that $C\ou{R^e}L$ is a finitely generated projective
  left $\ol S$-module. Since $D$ is a finitely generated projective
  left $R^e$-module, we can conclude that $\wtilde L=C\ou{R^e}L\ou{R^e}D$
  is a finitely generated projective left $\ol S$-module.

  However, we do not know in this situation whether 
  $\wtilde L^\vee$ and $\widetilde{L^\vee}$ are isomorphic 
  $\tS$-bialgebras. 

  Recall that the left dual $L^\vee$ is finitely generated projective
  as left $R$-module. For $\tR$-bialgebras $H$ such that $_RH$ is
  finitely generated projective, one can define a right dual
  $\tR$-bialgebra $^\vee H$ in such a way that $^\vee(L^\vee)\cong L$.

  Now let $\hat L:={^\vee\widetilde{L^\vee}}$ be the right
  dual $\tS$-bialgebra of the $\tS$-bialgebra obtained by 
  \sqm\ base change from the right dual of $L$ (note that
  $\widetilde{L^\vee}$ is a finitely generated projective left
  $S$-module by reasoning similar to that used for $_{\ol S}\wtilde L$).
  Then we have equivalences of monoidal categories
  $$\LComod L\cong\LMod{L^\vee}\cong\LMod{\widetilde{L^\vee}}
    \cong\LComod{\hat L}.$$
  
  If our \sqm\ equivalence comes from a Morita equivalence, then 
  $\hat L\cong\wtilde L$. Otherwise, it seems that we have another
  version of \sqm\ base change, suitable for comodules instead
  of modules. 

  More conjecturally, such a dual version of \sqm\ base
  change should also be possible if $L$ is not assumed to be
  finitely generated projective as left $\ol R$-module.
\end{Rem}

\section{Canonical Tannaka duality}\nmlabel{Section}{sec:ctd}
In this section we let $k$ be a field.
Let $\C$ be a semisimple $k$-linear tensor category with a finite
number of simple objects whose endomorphism rings are isomorphic
to $k$. 
Hayashi \cite{Hay:CTDFSTC,Hay:BIFA} has proved that $\C$ is 
equivalent to the category of 
modules over
a finite dimensional face algebra $F$.
The construction can of course be applied to $\LMod H$ where
$H$ is a split semisimple quasi-Hopf algebra, though it is not
so clear how $F$ is related to $H$.

In this section we will
describe a connection between the ``given'' $H$ and the 
``canonical'' $F$. This proceeds in two steps. First, one uses
a generalized smash product construction that produces a
$\times_H$-bialgebra $L$ isomorphic to $H\o H\o H^*$ as a vector
space, and with $\LMod L\cong\LMod H$ as monoidal categories.
In a second step, we use Morita base change to replace $H$ by
the Morita equivalent product of copies of the base field.
The result is a face algebra $\wtilde L$, and we shall show
that $\wtilde L\cong F$.

Let us first recall some elements of Hayashi's construction.
An important step is the construction of a monoidal functor
$\Omega_0\colon \C\rightarrow\BiMod R$, where $R=k^n$ and
$n$ is the number of isomorphism classes of simple objects
in $\C$. We will not go into details on the second important
step, which is the construction of a 
unique face algebra $F$ with $n$ face idempotents 
such that
$\Omega_0$ factors over a monoidal equivalence
$\Omega\colon\C\rightarrow\LMod F$. (In fact Hayashi uses modules
rather than comodules, which is of no importance since the face
algebra he constructs is finite dimensional.) 

Let $\Lambda$ be the set of isomorphism classes of simple 
objects in $\C$. 
For simplicity we let $R\cong k^n$ have the set $\Lambda$ as its
canonical basis of idempotents.

Hayashi's canonical functor $\Omega_0$ sends $X\in\C$
to the $R$-bimodule $\Omega_0(X)$ with
$\mu\Omega_0(X)\lambda=\Hom_{H-}(L_\mu,X\o L_\lambda)$,
where, compared with Hayashi's comvention, we have 
switched the sides in $R$-bimodules and replaced tensor product
in $\C$ by its opposite.
The monoidal functor structure
$\omega$ of $\Omega_0$
is the map
\begin{align*}
   \Omega_0(X)\ou R\Omega_0(Y)&\xrightarrow{\omega}\Omega_0(X\o Y)\\
   \mu\Omega_0(X)\rho\o\rho\Omega_0(Y)\lambda&\rightarrow\Omega_0(X\o Y)\\
    f\o g&\mapsto (L_\mu\xrightarrow{f}X\o L_\rho\xrightarrow{X\o g}X\o(Y\o L_\lambda)\xrightarrow{\Phi\inv}(X\o Y)\o L_\lambda),
\end{align*}
where $\mu,\rho,\lambda\in\Lambda$, and $\Phi$ denotes the associator
isomorphism in the category $\C$.

Now let $H$ be a split semisimple quasi-Hopf algebra. We will apply
Hayashi's constructions to $\C=\LMod H$, and investigate the relation
of $F$ to $H$.

The first step is a construction suggested by Hausser and Nill
(see \cite{HauNil:ITQHA}, Proposition 3.11 and the remarks following
the proof): 
They have defined a category $\HMod{}HHH$
of Hopf modules over $H$, 
which is monoidal in such a way that the underlying functor
$\mathcal U\colon\HMod{}HHH\rightarrow \BiMod H$ is a strict monoidal
functor.
They show $\HMod{}HHH$ to be equivalent 
as a monoidal category
to $\LMod H$ via a monoidal functor
$$(\mathcal R,\xi)\colon\LMod H\ni V\rightarrow V\o H\in\HMod{}HHH.$$

Now by translating
the coaction of $H$ on a Hopf module into an action of the dual $H^*$,
one can describe Hopf modules in $\HMod{}HHH$ equivalently as
modules over a certain generalized smash product
$L:=(H\o H^\op)\# H^*$. 
This kind of classification of Hopf modules by modules over an
algebra which is a product of several copies of $H$ and its dual
goes back to Cibils and Rosso \cite{CibRos:HBM}.
We refer the reader to \cite[Ex.4.12]{Sch:AMCGHSP} for
details on the construction of $L$.
Since the underlying functor
$$\LMod{(H\o H^\op)\#H^*}\cong\HMod{}HHH\rightarrow\BiMod H$$
is strictly monoidal, 
it follows from 
\cite[Thm.5.1]{Sch:BNRSTHB}
that $L$ has the structure of a
$\times_H$-bialgebra such that 
$\LMod{L}\cong\HMod{}HHH\cong\LMod H$ are equivalences
of monoidal categories. 

$H$ being split semisimple, it is 
Morita equivalent to a direct product of copies of $k$. We
claim that Hayashi's $F$ results from applying the appropriate
Morita base change to $L$.

To see this, we have to verify that the diagram
of monoidal functors
$$\xymatrix{\LMod H\ar[r]^-{\mathcal R}\ar[d]^{\Omega_0}&
            \HMod{}HHH\ar[d]^{\mathcal U}\\
            \BiMod R&\BiMod H\ar[l]^{\F}}$$
commutes up to isomorphism of monoidal functors. Here $\F$ denotes
the monoidal functor given by the Morita equivalence between $R$
and $H$.

For $\lambda\in\Lambda$, 
now the set of simple modules in $\LMod H=\C$,
fix a minimal idempotent
$e_\lambda\in H$ such that $L_\lambda:=He_\lambda\in\Lambda$.

The functor $\mathcal {UR}$ maps $X\in\LMod H$ to 
$\mathcal {UR}(X)={_\cdot X}\o{_\cdot H_\cdot}\in\BiMod H$; here
the dots indicate that $X\o H$ is equipped with the diagonal left
$H$-module structure and the right module structure induced by that
of the right tensor factor. The functor $\mathcal F$ maps 
$M\in\BiMod H$ to the $R$-bimodule defined by
$\mu\mathcal F(M)\lambda=e_\mu Me_\lambda$, so that
$\mathcal {FUR}(X)$ satisfies
\begin{multline*}
  \mu\mathcal{FUR}(X)\lambda=
   e_\mu\mathcal{UR}(X)e_\lambda=
   e_\mu(X\o H)e_\lambda
   \cong\Hom_{H-}(He_\mu,X\o He_\lambda)
   \\=\Hom_{H-}(L_\mu,X\o L_\lambda)
   =\mu\Omega_0(X)\lambda.
\end{multline*}
Note that the isomorphism $\psi\colon\Omega_0(X)\cong \mathcal F(X\o H)$ 
we have found maps $f\in\Hom_{H-}(L_\mu,X\o L_\lambda)$ to
$\psi(f)=f(e_\mu)\in e_\mu(X\o H)e_\lambda$.

We still have to show that $\Omega_0\cong\mathcal{FUR}$ as
{\em monoidal} functors.

The monoidal functor structure $\xi$ of $\mathcal R$
is the isomorphism
$$(X\o H)\ou H(Y\o H)=X\o(Y\o H)\xrightarrow{\Phi\inv}(X\o Y)\o H$$
in which the first equality is the canonical identification.

For $M,N\in\BiMod H$ we can 
identify $\F(M\ou HN)\subset M\ou HN$ with $\F(M)\ou R\F(N)$, which
makes $\F$ a strict monoidal functor. In particular, the monoidal
functor structure of $\mathcal{FUR}$ is the restriction of that of
$\mathcal R$; we shall denote it by $\xi$ again.

We need to show that the diagrams
$$\xymatrix{\Omega_0(X)\ou R\Omega_0(Y)\ar[r]^-{\omega}\ar[d]^{\psi\o\psi}
            &\Omega_0(X\o Y)\ar[d]^\psi\\
            \mathcal {FUR}(X)\ou R\mathcal {FUR}(Y)\ar[r]^-{\xi}
            &\mathcal {FUR}(X\o Y)}$$
commute for $X,Y\in\LMod H$. Let $f\in\Hom_H(L_\mu,X\o L_\rho)$ 
and $g\in\Hom_H(L_\rho,Y\o L_\lambda)$. Then 
$\psi\omega(f\o g)=\omega(f\o g)(e_\mu)=\Phi\inv(X\o g)f(e_\mu)$.
On the other hand, write
$f(e_\mu)=\sum x_i\o h_i$
with $x_i\in X$ and $h_i\in L_\rho$. 
Then we have
\begin{multline*}
 \xi(\psi\o\psi)(f\o g)
=\xi(f(e_\mu)\o g(e_\rho))=\Phi\inv(\sum x_i\o h_ig(e_\rho))
\\=\Phi\inv(\sum x_i\o g(h_ie_\rho))=\Phi\inv(\sum x_i\o g(h_i))
=\Phi\inv(X\o g)f(e_\mu).
\end{multline*}

\section{An example from Subfactor theory}\nmlabel{Section}{sec:II1}
Nikshych and \Vainerman\ have shown how to associate a weak Hopf 
algebra to a subfactor of finite depth 
of a von Neumann algebra factor  
\cite{NikVai:CD2SII1F,NikVai:GCII1FQG}.
The case of a subfactor $N\subset M$ of a type $\II_1$ factor of index
$\beta=4\cos^2\frac\pi{n+3}$ with $n\geq 2$ is treated in more detail
in \cite{NikVai:CD2SII1F,NikVai:FQGA}. 
The associated weak Hopf algebra
can be described as follows
(we summarize the beginning of
\cite[sec.2.7]{NikVai:FQGA}): Let $\mathcal A_{\beta,k}$ be the 
Temperley-Lieb algebra as in \cite{GooHarJon:CGTA}, that is, the
unital algebra freely generated by idempotents $e_1,\dots,e_{k-1}$
subject to the relations
$\beta e_ie_je_i=e_i$ for $|i-j|=1$ and $e_ie_j=e_je_i$ for 
$|i-j|\geq 2$. Then $\mathcal A_{\beta,k}$ is semisimple for 
$k\leq n+1$ by the choice of $\beta$ 
(cf.\ \cite[\S 2.8]{GooHarJon:CGTA}). Define
$A_{1,k}$ by $A_{1,k}=\mathcal A_{\beta,k+1}$ if $k\leq n+1$, 
and let $A_{1,k+1}$ be obtained by applying the Jones basic 
construction to the inclusion $A_{1,k-1}\subset A_{1,k}$ for
$k\geq n+1$. Thus $H:=A_{1,2n-1}$ is generated by idempotents
$e_1,\dots,e_{2n-1}$, and contains $A_{1,n-1}$, generated by
$e_1,\dots,e_{n-1}$, and $A_{n+1,2n-1}$, generated by
$e_{n+1},\dots,e_{2n-1}$, as subalgebras. Nikshych and \Vainerman\
describe a weak Hopf algebra structure of
$H$ with target counital subalgebra $H_t=A_{1,n-1}$ and
$H_s=A_{n+1,2n-1}$. For $n=2$, $H_t\cong \mathbb C\oplus \mathbb C$
is commutative, and $H\cong M_2(\mathbb C)\oplus M_3(\mathbb C)$
is a face algebra of dimension $13$. 

We shall examine the case
$n=3$. Here the Bratteli diagram of the inclusions
$A_{1,2}\subset\dots\subset A_{1,5}$ is obtained from that of
the inclusion $A_{1,2}\subset A_{1,3}$ (which is found in 
\cite{GooHarJon:CGTA}) by applying Jones' basic construction
twice:
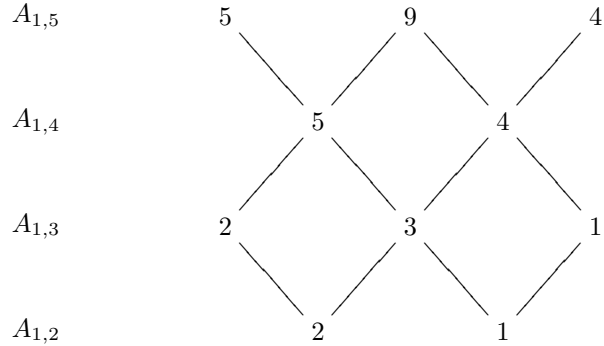
\newcommand\str{\ar@{-}}
$$\xymatrix{%
A_{1,5}&&   5\str[dr]&&
            9\ar@{-}[dl]\ar@{-}[dr]&&
            4\ar@{-}[dl]\\
A_{1,4}&&           &5\ar@{-}[dl]\ar@{-}[dr]&&
            4\ar@{-}[dl]\ar@{-}[dr]\\
A_{1,3}&&            2\str[dr]&&3\str[dl]\str[dr]&&1\str[dl]\\
A_{1,2}&& &2&&1}$$

We see that $H=A_{1,5}$ has dimension
$5^2+9^2+4^2=122$, and its counital subalgebras are
isomorphic to $\mathbb C\oplus M_2(\mathbb C)$, of dimension $5$.
We will apply Morita base change to this example, reducing the 
counital subalgebra to $\mathbb C\oplus \mathbb C$. We will not
derive an explicit description of the resulting face algebra, but
will be content with determining its algebra structure, hence 
its dimension.

The Bratteli diagram for the inclusion $H_t\subset H$ is
obtained by composing the stories of the Bratteli diagram above:
$$
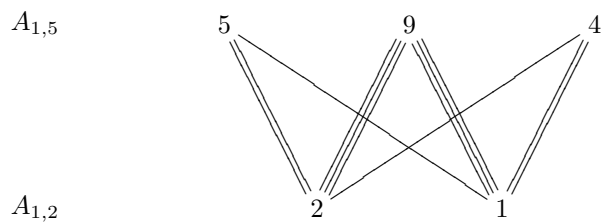
\xymatrix{
A_{1,5}&&5\ar@{=}[ddr]\ar@{-}[ddrrr]&&9\ar@3{-}[ddl]\ar@3{-}[ddr]&&4\ar@{-}[ddlll]\ar@{=}[ddl]\\\\
A_{1,2}&&            &2&&1}$$
To apply Morita base change, we need the Bratteli diagram for the
map $H_s\otimes H_t\rightarrow H$.
Of course $H_s\otimes H_t\cong M_4(\mathbb C)\oplus M_2(\bC)\oplus M_2(\bC)\oplus\bC$,
and the antipode of $H$, which exchanges $H_s$ and $H_t$, will 
switch the two copies of $M_2(\bC)$. 
The lower story of the following tower is the inclusion $H_t\subset H_t\o H_s$:

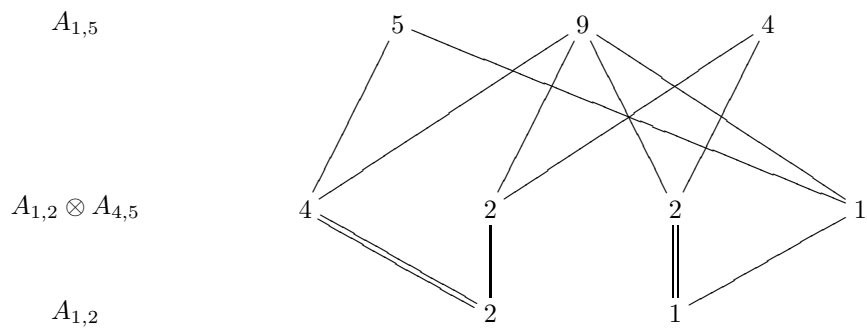
$$\xymatrix{
A_{1,5}&&   &5\ar@{-}[ddl]\ar@{-}[ddrrrrr]
            &&9\ar@{-}[ddlll]\ar@{-}[ddl]\ar@{-}[ddr]\ar@{-}[ddrrr]
            &&4\ar@{-}[ddl]\ar@{-}[ddlll]\\\\
A_{1,2}\o A_{4,5}&&            4\ar@{=}[drr]&&2\ar@{-}[d]&&2\ar@{=}[d]&&1\ar@{-}[dll]\\
A_{1,2}&&            &&2&&1}$$
We claim that the upper story is the Bratteli diagram
for the map $H_t\o H_s\rightarrow H$. It will help to know that, 
since the antipode switches the two vertices labelled $2$ on the
middle floor, any top floor vertex has the same number of edges
to each of these vertices.
In particular, there can be no edge from the $5$ 
on the top floor to a $2$ on the middle floor, since
there is only one $1$ on the bottom floor linked to $5$. Also, there
can be no more than one link from $5$ to the $1$ on the middle floor,
since there should be only one link to the $1$ on the bottom floor.
This makes the two links from $5$ to the middle floor inevitable
as shown. There should be three links from the top $9$ to the
bottom $2$. These cannot be accounted for by three links to the
left $2$ on the middle floor, since that would also entail 
three links to the right $2$, and hence six links to the bottom
$1$. So there has to be one link to the $4$ and one to the left $2$,
hence also the right $2$ on the middle floor, which makes the
one link to the $1$ on the middle floor also inevitable. Finally the one
link from the top $4$ to the bottom $2$ can only be accounted for
by a single link from the top $4$ to the left $2$ on the middle
floor, hence there also has to be one to the right $2$, and there
is no room for more.

Now we apply Morita base change 
to pass from $H$ with counital subalgebra 
$H_t\cong \bC\oplus M_2(\bC)$ to $\wtilde H$ with counital 
subalgebra $\wtilde H_t\cong\bC\oplus\bC$. The Bratteli diagram for the
inclusion $\wtilde H_t\oplus \wtilde H_s\subset \wtilde H$ is the same as for
$H_t\oplus H_s\subset H$, but with all ranks on the lower
floor replaced by $1$:
$$
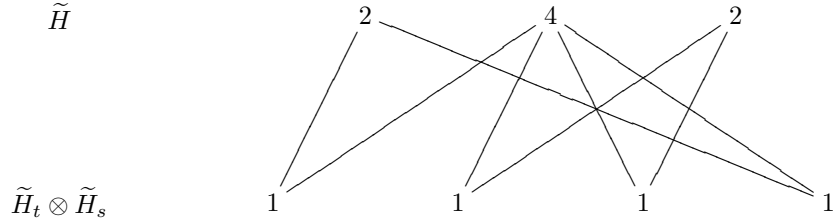
\xymatrix{
\wtilde H&&   &2\ar@{-}[ddl]\ar@{-}[ddrrrrr]
            &&4\ar@{-}[ddlll]\ar@{-}[ddl]\ar@{-}[ddr]\ar@{-}[ddrrr]
            &&2\ar@{-}[ddl]\ar@{-}[ddlll]\\\\
\wtilde H_t\o \wtilde H_s&&            1&&1&&1&&1}$$
The resulting ranks on the upper floor show that 
$\dim\wtilde H=2^2+4^2+2^2=24$. (Remember that $\dim H=122$.)
\begin{Rem}
  By Hayashi's canonical Tannaka duality, there is a face algebra
  $F$ with $\LMod{F}\cong\LMod H$ as monoidal categories, 
  where $F$ has three face idempotents (since $H$ has three
  isomorphism classes of irreducible modules) whereas our $\wtilde H$
  has two. Hayashi has also described another procedure to associate
  a face algebra to any subfactor of a $\II_1$ factor 
  \cite{Hay:GQGII1S},
  which will, however, also yield a face algebra that has as many
  faces as isomorphism classes of irreducible modules. 
  By contrast, applying Morita base change to the weak Hopf algebra
  $A_{1,2n-1}$
  of Nikshych and \Vainerman\ will yield a face algebra with one 
  face less than isomorphism classes of irreducible modules whenever
  $n$ is odd.
\end{Rem}

\end{document}